\documentclass{amsart}
\usepackage{amssymb, mathtools, fullpage, color, cite, float}
\usepackage[colorinlistoftodos]{todonotes}
\usepackage[colorlinks=true, pdfstartview=FitV,linkcolor=blue,citecolor=blue,urlcolor=blue]{hyperref}
\usepackage{bbm}

\theoremstyle{definition}
\newtheorem {theorem}{Theorem}[section]

\newtheorem{remark}[theorem]{Remark}
\newtheorem{example}[theorem]{Example}

\newenvironment{red}{\relax\color{red}}{\relax}
\newenvironment{blue}{\relax\color{blue}}{\hspace*{.5ex}\relax}

\newcommand{\ber}{\begin{red}}
\newcommand{\er}{\end{red}}
\newcommand{\beb}{\begin{blue}}
\newcommand{\eb}{\end{blue}}

\numberwithin{equation}{section}

\begin{document}

\title[Computer vision and converse theorem]{Computer vision and converse theorems}

\date{\today}

\author[Y.-H. He]{Yang-Hui He}
\address{London Institute for Mathematical Sciences, Royal Institution, London W1S 4BS, UK  \hfill \break \indent Merton College, Oxford, OX14JD, UK}
\email{\href{mailto:hey@maths.ox.ac.uk}{hey@maths.ox.ac.uk}}

\author[K.-H. Lee]{Kyu-Hwan Lee}
\address{Department of Mathematics, University of Connecticut, Storrs, CT 06269, USA  \hfill \break \indent Korea Institute for Advanced      Study, Seoul 02455, Republic of Korea}
\email{\href{mailto:khlee@math.uconn.edu}{khlee@math.uconn.edu}}

\author[T. Oliver]{Thomas Oliver}
\address{University of Westminster, London, UK}
\email{\href{mailto:T.Oliver@westminster.ac.uk}{T.Oliver@westminster.ac.uk}}

\author[Y. Qi]{Yidi Qi}
\address{Department of Physics, Northeastern University, Boston, MA, USA \hfill \break \indent NSF Institute for Artificial Intelligence and Fundamental Interactions, Cambridge, MA, USA}
\email{\href{mailto:y.qi@northeastern.edu}{y.qi@northeastern.edu}}

\begin{abstract}
Random matrices provide a well-established statistical model for a range of arithmetic phenomena.
In this paper, we investigate the extent to which one- and two-dimensional convolutional neural networks (CNNs) can distinguish between arithmetic data arising from elliptic curves with conductor in a fixed interval and random matrix data drawn from the same Sato--Tate distribution.
Inspired by converse theorems in the Langlands program, we represent each elliptic curve together with its twists as a vector field and, subsequently, encode that vector field as a digital image.
We observe that a two-dimensional CNN trained on this image data is better able to separate conductor families from random matrix data than a one-dimensional CNN trained on vectors of Frobenius traces without twisting data.
We also observe that the same two-dimensional architecture can predict the analytic rank of an elliptic curve, and it does so by factoring through the (untwisted) Frobenius traces.
\end{abstract}

\maketitle

\section{Introduction}\label{s.intro}
With the rapid development of artificial intelligence over the past decade, we have entered an era of AI-assisted mathematical research \cite{He:2017aed,douglas2022machine,Gukov:2024buj,he2024can,he2024ai,douglas2025mathematical}. 
Techniques that were once unfamiliar to pure mathematicians are now increasingly prominent, particularly in number theory, where data-driven and statistical approaches have proved especially fruitful \cite{Alessandretti:2019jbs,He:2020kzg,He:2020nf,He:2020tkg, He:2020ahlos}. 
In this paper, we study a fundamental problem in arithmetic from a novel perspective using convolutional neural networks (CNNs).

To each elliptic curve over $\mathbb{Q}$ one may associate a sequence of Frobenius traces. 
In this paper, we restrict attention to the generic case of elliptic curves without complex multiplication. 
Statistically, these traces are modeled by the traces of random matrices, generically from $\mathrm{SU}(2)$. 
It should not be easy to distinguish arithmetic data from random matrix data with the same distribution.
Nonetheless, we achieve the distinction by supplementing statistics with an insight from the representation theory of automorphic forms.

By modularity, the $L$-function of an elliptic curve over $\mathbb{Q}$ coincides with that of a modular form. 
Weil's converse theorem (and its variants) provides a theoretical mechanism for distinguishing the coefficients of modular $L$-functions from arbitrary complex sequences \cite{weil}. 
The key idea is to consider not only the original sequence of coefficients, but also sufficiently many twisted sequences parameterized by Dirichlet characters.
This information is naturally presented as a two-dimensional array which may be encoded as a digital image.
Our central observation is that image recognition techniques can exploit this structure to distinguish elliptic curves of small conductor from random matrix data. Based on this observation, the main message of this paper is that converse theorems can be implemented as instances of image classification.

CNNs are a class of machine learning models widely used in computer vision and pattern recognition \cite{Goodfellow-et-al-2016}. 
The possibility that problems in algebraic geometry and number theory might be recast in this framework was proposed in \cite{He:2017aed,He:2018jtw}, although issues of interpretability remain challenging. 
A particularly rich source of arithmetic data is the LMFDB \cite{lmfdb}, which has played a central role in several AI-driven studies of number theory \cite{Alessandretti:2019jbs,He:2020kzg,He:2020nf,He:2020tkg, He:2020ahlos}. These efforts culminated in the discovery of murmurations \cite{He:2022pqn}, which are now an active area of research \cite{bober2023murmurations,zubrilina2024convergence,booker2024murmurations,shi2024murmurations,cowan2024murmurations,babei2025learning,lee2025murmurations,Bieri:2025tyk,cowan2025mean,bujanovic2025improving,lee2025machines,sawin2025murmurations,kuan2025murmurations,lowry2025murmurations,martin2025variations,shi2025murmurations,martin2025distribution,wang2025murmurations}.

We conclude this introduction with a brief outline of the paper. 
In Section~\ref{s.prelim}, we introduce the necessary mathematical background regarding elliptic curves. 
Section~\ref{s.ima} describes the construction of our dataset of two-dimensional arrays and their encoding as digital images. 
In Section~\ref{s.cnn}, we present the architecture of our CNNs and discuss the experimental results. 
In particular, we demonstrate that a one-dimensional CNN is less able to distinguish arithmetic data from random matrices than an analogous two-dimensional CNN incorporating twist. 
Subsequently, we observe that the same two-dimensional CNN can distinguish unseen elliptic curves of larger conductor from random matrix data, and that it can accurately predict the rank of an elliptic curve.
Throughout, we exhibit saliency maps in order to visualise and interpret how each prime/twist contributes to the predictions of the neural networks.

\subsection*{Acknowledgements}
YHH thanks the Leverhulme Trust Research Project (Grant No. RPG-2022-
145) and STFC grant ST/J00037X/4. 
YDQ was supported by the NSF grant PHY-2019786 (the NSF AI Institute for Artificial Intelligence and Fundamental Interactions).

\section{Preliminaries}\label{s.prelim}
We begin with the requisite number theory used in the sequel.
\subsection{Elliptic curves over $\mathbb{Q}$}
An elliptic curve $E$ over $\mathbb{Q}$ has a Weierstrass model of the form: \begin{equation}\label{weier5}
E:~y^2+w_1xy+w_3y=x^3+w_2x^2+w_4x+w_6.
\end{equation}
Each model for $E$ has a discriminant $\Delta$, and a minimal Weierstrass equation for $E$ is one such that $|\Delta|$ is minimal among all models of $E$.
We say that a prime number $p\in\mathbb{Z}$ is a good (resp. bad) prime for $E$ if the reduction mod $p$ of a minimal Weierstrass equation for $E$ is non-singular (resp. singular).
The conductor $\mathcal{N}(E)$ is the product over all the primes of bad reduction, each to an appropriate power that indicates the reduction type \cite{silverman2013advanced}.

The set $E(\mathbb{Q})$ of rational points on an elliptic curve $E$ admits the structure of a finitely generated abelian group. 
We define the rank $r(E)$ to be the rank of this group.
For each prime $p$ and integer $k>0$, the number of points $\#E(\mathbb{F}_{p^k})$ of $E$ over the finite field $\mathbb{F}_{p^k}$ is finite. 
For each prime $p$, we define the local Hasse--Weil zeta function:
\begin{equation}
    Z_p(E, s) :=
    \exp \left( \sum\limits_{k=1} \frac{\#E(\mathbb{F}_{p^k})}{k} p^{-s k}
    \right) \ .
\end{equation}
The local Hasse--Weil zeta function is known to be a rational function in $p^{-s}$, whose logarithmic derivative is a generating function for the sequence $\#E(\mathbb{F}_{p^k})$ \cite{silverman2013advanced}.
Taking the product over all primes gives the Hasse--Weil zeta function
\begin{equation}
    Z(E,s) = \prod\limits_p Z_p(E, s)
    = \frac{\zeta(s) \zeta(s-1)}{L(E, s)}, \ \ \mathrm{Re}(s)\gg0.
\end{equation}
The denominator is called the elliptic $L$-function, and takes the form
\begin{equation}\label{eq.elf}
L(E,s)=\prod_{p:\text{ good prime}}\left(1-a_p(E)p^{-s}+p^{1-2s}\right)^{-1}\prod_{p:\text{ bad prime}}\left(1-a_p(E)p^{-s}\right)^{-1},
\ \ \mathrm{Re}(s)\gg0.
\end{equation}
In the above, the coefficients $a_p$ are the so-called {\em Frobenius traces}
where, for a good prime $p$, we have 
\begin{equation}
    a_p(E)=p+1-\#E(\mathbb{F}_p) 
\end{equation}
and, for a bad prime $p$, we have $a_p(E)\in\{-1,0,1\}$ depending on the reduction type \cite{silverman2013advanced}.
Two elliptic curves $E_1,E_2$ are isogenous if and only if the sequences $(a_p(E_1))_p$ and $(a_p(E_2))_p$ agree almost always.
There are some key properties of $a_p(E)$ relevant to the present study. 
Firstly, for a prime $p$, the Hasse bound dictates that
\begin{equation}\label{hasse}
    |a_p(E)| \leq 2 \sqrt{p} \ . 
\end{equation}
Given equation \eqref{hasse}, it is natural to 
define $x_p := a_p/2\sqrt{p} \in [-1,1]$.
If $E$ has no complex multiplication (CM), then the sequence \( (x_p)_{p \nmid \mathcal{N}(E)} \) 
is equidistributed in the interval \([-1,1]\) with respect to the \emph{Sato–Tate measure}
\begin{equation}\label{eq.STnCM}
\mu(dx) = \frac{2}{\pi} \sqrt{1 - x^2}\, dx.
\end{equation}
There is an analogous distribution in the case that $E$ has CM \cite{serre2016lectures}.

\subsection{Converse theorems}
For $M\in\mathbb{Z}_{>0}$, a Dirichlet character mod $M$ is a completely multiplicative function $\chi:\mathbb{Z}\rightarrow\mathbb{C}$ that is periodic with period $M$ and satisfies $\chi(a)=0$ if and only if $\mathrm{gcd}(a,M)>1$. 
The conductor of a Dirichlet character $\chi$ modulo $N$ is the least positive integer $q$ dividing $M$ for which $\chi(n+kq)=\chi(n)$ for all $n$ and $n+kq$ coprime to $M$.
We say that $\chi$ is primitive if its modulus is equal to its conductor.
For example, the principal character $\chi_0$ mod $M$, given by $\chi_0(a)=1$ (resp. $\chi_0(a)=0$) if $\mathrm{gcd}(a,M)=1$ (resp. $\mathrm{gcd}(a,M)>1$), is not primitive.
If $M=q$ is prime, then all non-principal characters mod $q$ are primitive, and the set of characters mod $q$ form a cyclic group with identity $\chi_0$.

If $E$ is an elliptic curve over $\mathbb{Q}$ with conductor $\mathcal N$, then there exists a weight $2$ modular form $f$ on $\Gamma_0(\mathcal N)$ such that $L(f,s)=L(E,s)$.
Subsequently, it follows that $L(E,s)$ admits analytic continuation to $\mathbb{C}$ and satisfies a well-known functional equation \cite{bump1998automorphic}.
Conversely, one can ask what properties must hold for a sequence $(a_n)_{n=1}^{\infty}$ to determine a modular $L$-function $L(s)=\sum_na_nn^{-s}$.
The prototypical answer is given by Weil's converse theorem \cite{weil}, which involves twisted $L$-functions $L(s,\chi)=\sum_{n=1}^{\infty}\chi(n)a_nn^{-s}$,
in which $\chi$ is a primitive Dirichlet character.
Although we will not need the technical details of Weil's converse theorem we are inspired by the underlying principle, namely, if these twisted $L$-functions satisfy analytic continuation, functional equations, and suitable growth conditions for sufficiently many $\chi$, then the sequence $(a_n)_{n=1}^{\infty}$ arises from a modular form.
In this paper, we will use a framework inspired by Weil's converse theorem to distinguish elliptic $L$-functions with small conductor from random matrix data with the same distribution.

\section{Arithmetic families as digital images}\label{s.ima}
In this Section, we will explain how our datasets were generated.
The code required to reproduce the data is available at \cite{repo}.

\subsection{Conductor families}\label{s.confuctorfamily}
By a \textsl{conductor family}, we mean arithmetic objects with similar conductors, possibly sharing other common features (for example, rank).
Conductor families were, for instance, very fruitful in the data-driven discovery of murmurations \cite{He:2022pqn}.

For $S\subset\mathbb{Z}_{\geq0}$ and an interval $I\subset\mathbb{R}_{>0}$ we may consider the set $\mathcal{E}_S(I)$ of isogeny classes with rank $r\in S$ and conductor in $I$.
For an interval $J\subset\mathbb{R}_{>0}$, we define, from the Frobenius traces $a_p(E)$: 
\begin{equation}\label{eq.cmatrix}
C^{\text{matrix}}(S,I,J)=\left(\frac12-\frac{a_p(E)}{4\sqrt{p}}\right)_{E\in\mathcal{E}_S(I),p\in J},
\end{equation}
in which rows are indexed by curves and columns are indexed by primes.
Using the Hasse bound in \eqref{hasse}, we see that the coefficients of $C^{\text{matrix}}(S,I,J)$ are in $[0,1]$. 
By rescaling and rounding, we may interpret each matrix element as a pixel in an 8-bit greyscale image, and its value as determining the corresponding shade of grey. 
Indeed, we can define
\begin{equation}\label{eq.cfimage}
    C^{\text{image}}(S,I,J)=
    \lfloor (2^8 - 1)C^{\text{matrix}}(S,I,J) \rfloor
=
\left(\left\lfloor
    (2^8 - 1)\left(\frac12-\frac{a_p(E)}{4\sqrt{p}}\right)\right\rfloor\right)_{E\in\mathcal{E}_S(I),p\in J}
\ ,
\end{equation}
where $\lfloor x \rfloor$ is the greatest integer which is $\leq x$. 
The significance of $2^8-1=255$ in equation~\eqref{eq.cfimage} is that each $8$-bit binary string in $\{0,1\}^8$ corresponds to an integer in the set $\{0,\dots,255\}$.
\begin{example}
In Figure~\ref{fig.rank}, we visualise equation~\eqref{eq.cfimage} as an $8$-bit greyscale image
for $\emptyset\neq S\subseteq\{0,1\}$ and $I=J=[1,1000]$.
\begin{figure}[htbp]
    \centering
    \includegraphics[width=0.9\textwidth]{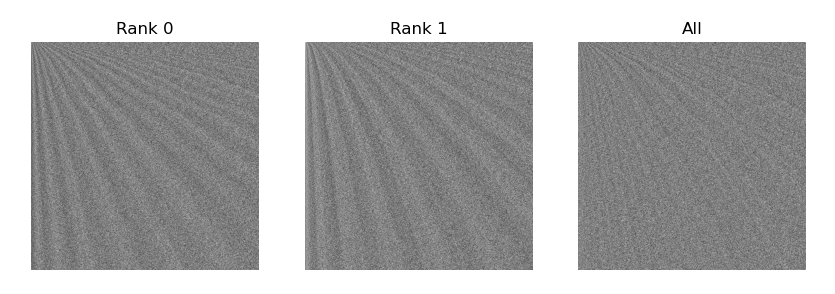}
    \caption{$8$-bit greyscale visualisation of matrix $(a_p(E))_{p,E}$, where $p$ varies over primes $<1000$ and $E$ varies over elliptic curves over $\mathbb{Q}$ with rank in set $\{0\}$, resp. $\{1\}$, $\{0,1\}$ and conductor $\mathcal{N}(E)\leq1000$. Image encoding is given by equation~\eqref{eq.cfimage}.}
     \label{fig.rank}
\end{figure}
\end{example}
\begin{remark}
If $E_1$ and $E_2$ are elliptic curves over $\mathbb{Q}$ such that $|a_q(E_1)-a_q(E_2)|=1$ for some sufficiently large prime $q$, then we have
\[\left\lvert\frac{\left(2^8-1\right)a_q(E_1)}{4\sqrt{q}}-\frac{\left(2^8-1\right)a_q(E_2)}{4\sqrt{q}}\right\rvert<1.\]
In other words, a small enough change in an element of $C^{\mathrm{matrix}}(S,I,J)$ is lost to rounding in the image encoding $C^{\mathrm{image}}(S,I,J)$. 
Similar loss due to image encoding is present in the next section, so, when necessary, we will work with raw numerical values such as those in equation~\eqref{eq.cmatrix}.
\end{remark}

\subsection{Twist families}\label{s.twist}
In this Section, we will encode each individual elliptic curve and finitely many twists as a vector field or image, in which the rows are indexed by primes and the columns are indexed by Dirichlet characters.
In this way, we will embed twisted $L$-functions into our datapoints.
In particular, we will visualise character twist families by coloured images.
More precisely, we will use true colour 24-bit images. 

Given an elliptic curve $E$, a primitive Dirichlet character $\chi$, and a prime $p$, we write
\begin{equation}\label{eq.reim}
\frac{a_{p}(E)\chi(p)}{2\sqrt{p}}=X_{p,\chi}+iY_{p,\chi}, \qquad X_{p,\chi}, Y_{p,\chi} \in \mathbb R.
\end{equation}
In other words $X_{p,\chi}$ (resp. $Y_{p,\chi}$) is the real (resp. imaginary) part of $a_p(E)\chi(p)/2\sqrt{p}$.
Note that:
\[\left\lvert X_{p,\chi}+iY_{p,\chi}\right\rvert=\sqrt{X^2_{p,\chi}+Y^2_{p,\chi}}<1,\]
and so we have $X_{i,j},Y_{i,j}\in[-1,1]$.
Therefore, the following values are in $[0,1]$:
\[R_{p,\chi}=\frac12-\frac12X_{p,\chi}, \ \ B_{p,\chi}=\frac12-\frac12Y_{p,\chi}.\]
Similarly to Section~\ref{s.confuctorfamily}, in order to make $8$-bit values, we need to rescale and round:
\begin{equation}\label{eq.RBtild}
\begin{split}
\widetilde{R}_{p,\chi}=\lfloor (2^8 - 1)R_{p,\chi}\rfloor=\left\lfloor (2^8 - 1)\left(\frac12-\mathrm{Re}\left(\frac{a_{p}(E)\chi(p)}{4\sqrt{p}}\right)\right)\right\rfloor, \\
\widetilde{B}_{p,\chi}=\lfloor (2^8 - 1)B_{p,\chi}\rfloor=\left\lfloor (2^8 - 1)\left(\frac12-\mathrm{Im}\left(\frac{a_{p}(E)\chi(p)}{4\sqrt{p}}\right)\right)\right\rfloor.
\end{split}
\end{equation}
Equation~\eqref{eq.RBtild} accounts for $8+8=16$ bits of our $24$ bit pixels, and determine the red (R) and blue (B) contributions to each pixel.
The remaining $8$ bits will be set to a constant value.
More precisely, we may construct a true colour image in which each pixel is determined by the three-dimensional vector:
\begin{equation}\label{eq.rgb}
(\widetilde{R}_{p,\chi},\widetilde{G}_{p,\chi},\widetilde{B}_{p,\chi})=\left(\lfloor (2^8-1) R_{p,\chi}\rfloor,2^7-1,\lfloor (2^8-1)B_{p,\chi}\rfloor\right).
\end{equation}
In particular, the green (G) contribution is taken to be $2^7-1$. 

Equivalently, our true colour image may be viewed as a vector field on the two-dimensional lattice spanned by $p$ and $\chi$.
By ordering $p$ and $\chi$, we may realise this lattice within $\mathbb{Z}^2$ and define the vector field $V:\mathbb{Z}^2\rightarrow\mathbb{R}^3$ given by:
\begin{equation}\label{eq.vftci}
V(p,\chi)=(\widetilde{R}_{p,\chi},\widetilde{G}_{p,\chi},\widetilde{B}_{p,\chi})
\end{equation}
We refer to the three scalar components of 
$V$ as the red, green, and blue channels, and to their values as the channel intensities.
The term \textsl{spatial information} refers to dependence of the channel intensity on the pixel coordinate $(p,\chi)\in\mathbb{Z}^2$, so that variation in the spatial information across the lattice encodes the geometric structure of the image.
In our case, only the red and blue channels carry spatial information.
\begin{remark}
We note that the construction above is applicable to any finite sequence $(z_p)_{p<b}\subset\mathbb{C}$ satisfying $|z_p|<1$ for all $p<b$.
Indeed, it suffices to replace the left-hand side of equation~\eqref{eq.reim} by $z_p$.
This will be relevant in Section~\ref{s.2dcnn}.
\end{remark}
\begin{example}\label{ex.tcim}
In Figure~\ref{fig.CT1}, we visualise three 24-bit true colour images associated to $y^2+y=x^3-x^2$ using the construction outlined above. 
In Figure~\ref{fig.CT1}, $p$ varies over the first $100$ (resp. $200$, $300$ primes) and $\chi$ varies over primitive Dirichlet characters of conductor $\leq24$ (resp. $\leq32$, $\leq 53$).
This yields an image of size $100\times100$ (resp. $200\times200$, $300\times300$).

\begin{figure}[htbp]
    \centering
    \includegraphics[width=0.3\textwidth]{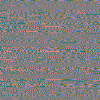}
    \includegraphics[width=0.3\textwidth]{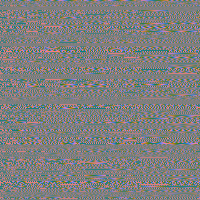}
    \includegraphics[width=0.3\textwidth]{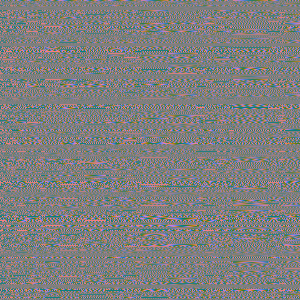}
    \caption{24-bit true colour images attached to  character twists of $y^2+y=x^3-x^2$ with resolutions (left) $100\times100$ (centre) $200\times200$ (right) $300\times300$. For details, refer to Example~\ref{ex.tcim}.}
     \label{fig.CT1}
\end{figure}
\end{example}

\section{Convolutional neural networks}\label{s.cnn}
In this Section, we will apply the techniques of computer vision to the vector fields attached to elliptic curves in Section~\ref{s.twist}.
The code required to reproduce the experiments of this paper is available at \cite{repo}.
For the general theory of CNNs, we refer the reader to \cite{Goodfellow-et-al-2016} or \cite{ImNet}.
A neural network is a finite composition of affine maps between finite-dimensional vector spaces and, typically non-linear, activation functions. 
An  $n$-dimensional CNN is a neural network in which some of the affine maps are convolution operators on functions on an $n$-dimensional lattice, typically followed by nonlinear activation functions and possibly other operations such as pooling (local aggregation) or normalisation (rescaling). 
The input for an $n$-dimensional CNN is a function $\Omega\rightarrow\mathbb{R}^c$, where $\Omega\subset\mathbb{Z}^n$ and $c$ is the number of channels.
Equivalently, the input is an $(n+1)$-dimensional tensor, $n$ spatial dimensions, and $1$ channel dimension of size $c$.
In this paper, the output of our CNNs will be a scalar computed by repeatedly reducing the spatial dimensions (resolution) of the input whilst simultaneously increasing the number of channels.
The scalar output is interpreted as the predicted category of the input.

In what follows, we report that a one-dimensional CNN is less effective at distinguishing elliptic curves of small conductor from random matrices with the same distribution than a two-dimensional CNN is at distinguishing twisted elliptic curves of small conductor from twisted random matrix data with the same distribution.

\subsection{One-dimensional CNN}\label{s.1dcnn}
Here, we consider a one-dimensional CNN with inputs determined by Frobenius traces and random matrix eigenvalues.
Let us define our dataset more precisely.

We work with the complete dataset $\mathcal{E}$ of 64,687 non-CM elliptic curves over $\mathbb{Q}$ with conductor $<10{,}000$. 
Code to generate the dataset $\mathcal{E}$ is available at \cite{repo}. 
We represent each $E\in\mathcal{E}$ as an $N$-dimensional vector
\begin{equation}\label{eq.vE}
v(E)=(a_2(E),a_3(E),\dots,a_{p_N}(E))\in\mathbb{R}^N \ ,
\end{equation}
where $p_n$ is the $n$-th prime number.
Equation~\eqref{eq.vE} is the same representation as used in \cite{He:2020tkg,He:2022pqn}.
For $N$, we will typically take values in the set $\{100,200,300\}$.
This yields a high-dimensional dataset $\mathcal{D}_1=\{v(E):E\in\mathcal{E}\}\subset\mathbb{R}^N$.
Next, we generate a dataset $\mathcal{D}_2$ of sequences $(\tilde{x}_p)_p\subset\mathbb{R}^N$ 
by drawing each coefficient independently from the Sato–Tate distribution as per equation~\eqref{eq.STnCM}.
More precisely, we use the following rejection sampling procedure.
First, fix a prime $p$. 
Sample $\theta_p$ uniformly from $[0,\pi]$ and $y$ uniformly from $[0,2/\pi]$.
If $y < \frac{2}{\pi}\sin^2\!\left(\theta_p\right)$, accept $\theta_p$; otherwise reject it and resample.
For an accepted $\theta_p$, define $\widetilde{x}_p=\cos(\theta_p)$ and $x_p=[2\widetilde{x}_p\sqrt{p}]$, where for $r\in\mathbb{R}$, $[r]$ denotes the integer obtained by rounding $r$.
For each prime $p$, we repeat this procedure to generate our dataset $\mathcal{D}_2$ consisting of $500{,}000$ vectors $(\tilde{x}_p)_p$.

The datasets $\mathcal{D}_1$ and $\mathcal{D}_2$ may be generated using code available at \cite{repo}. 
We note that the sets $\mathcal{D}_1$ and $\mathcal{D}_2$ are disjoint and define $\mathcal{D}=\mathcal{D}_1\sqcup\mathcal{D}_2$.
More philosophically, despite the fact that all datapoints in $\mathcal{D}_1$ and $\mathcal{D}_2$ adhere to the same distribution, only $\mathcal{D}_1$ is a meaningful conductor family.
In principle, one can try to associate a vector in $\mathcal{D}_2$ with an elliptic curve by searching for Weierstrass equations with the correct point counts modulo each prime and then using the Chinese remainder theorem to find a corresponding model over $\mathbb Z$. However, this process causes the Weierstrass coefficients to grow superexponentially and is not really feasible for constructing genuine curves. 

\begin{figure}[htbp]
    \centering
    \includegraphics[width=0.9\textwidth]{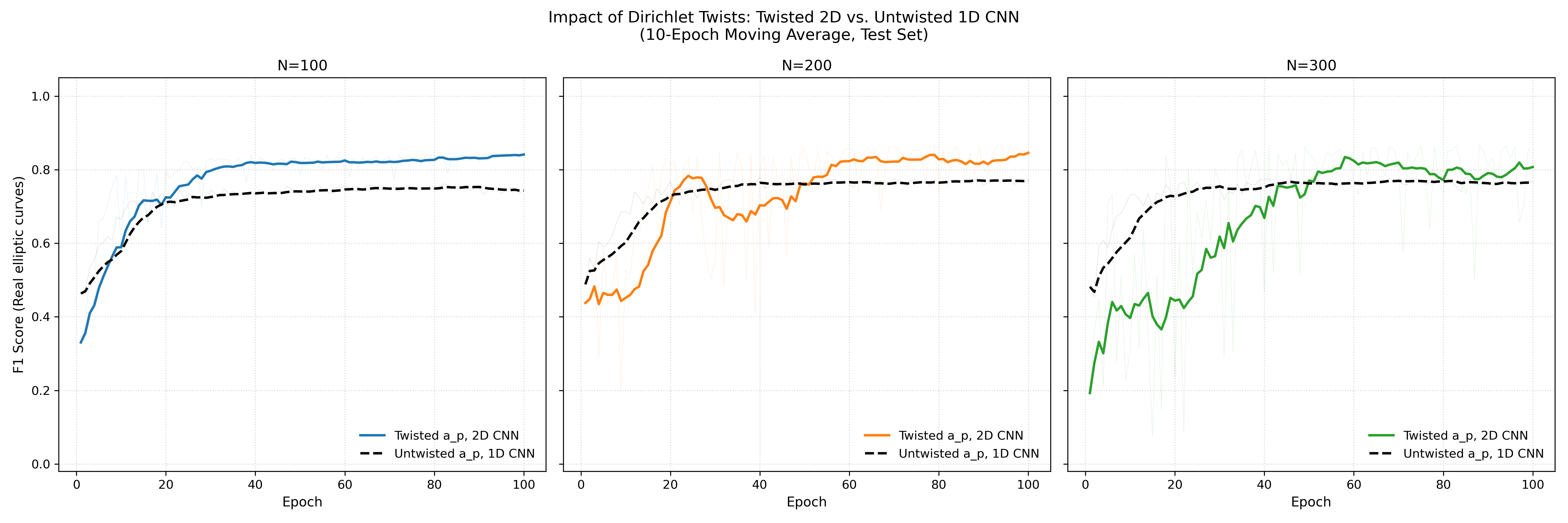}
    \caption{Evolution of F1 score relative to epoch for 1d CNN described in Section~\ref{s.1dcnn} and 2d CNN described in Section~\ref{s.2dcnn}.}
     \label{fig:f1_1d_2d_comparison}
\end{figure}

We set-up a simple binary classification problem: given $v\in \mathcal{D} \subset\mathbb{R}^N$, determine $i$ such that $v\in\mathcal{D}_i$.
To that end, we use a 1d CNN with architecture as in Table~\ref{tab:1dCNN} and apply an 80:20 training ratio. 
In words, after each convolutional layer, we apply batch normalisation, ReLU activation, and then max pooling to halve the spatial dimensions. After the final convolutional block, global average pooling reduces each feature map to a single scalar.
This is followed by a fully connected layer mapping
$512\mapsto 256$ features with batch normalisation, $\mathrm{ReLU}$ activation, and dropout (probability 0.5) for regularisation.
A final linear layer produces a single logit for binary classification.
The network is trained using the Adam optimizer with a learning rate of $0.001$.
Given the class imbalance between $\mathcal{D}_1$ and $\mathcal{D}_2$, we utilize a weighted binary cross entropy with logits loss function, assigning a higher weight of 3.0 to the minority class (the genuine elliptic curves).
For the same reason, we adopt the F1 score as our primary evaluation metric, which balances precision and recall. 
We note that standard accuracy would be misleadingly high for our dataset, which is dominated by $\mathcal{D}_2$. 
In Figure~\ref{fig:f1_1d_2d_comparison} (dotted line), we plot the evolution of the F1 score relative to the epoch for different choices of $N$.

\begin{table}[h]
\centering
\begin{tabular}{|c|c|c|}
\hline
Layer & Operation & Output shape \\
\hline
Input & -- & $(1, N)$ \\
Conv1 & $3$ conv, 64 ch & \\
Pool1 & $2$ max pool & $(64, [N/2])$ \\
Conv2 & $3$ conv, 128 ch & \\
Pool2 & $2$ max pool & $(128, [N/4])$ \\
Conv3 & $3$ conv, 256 ch &  \\
Pool3 & $2$ max pool & $(256, [N/8])$ \\
Conv4 & $3$ conv, 512 ch & \\
Pool4 & $2$ max pool & $(512, [N/16])$ \\
Conv5 & $3$ conv, 512 ch & \\
Pool5 & $2$ max pool & $(512, [N/32])$ \\
GAP & Global average pooling & $(512)$ \\
FC1 & Fully connected & $(256)$ \\
FC2 & Fully connected & $(1)$ \\
\hline
\end{tabular}
\caption{1d CNN architecture used for classification problem described in Section~\ref{s.1dcnn}. Batch normalisation and $\mathrm{ReLU}$ activations follow each convolution.}
\label{tab:1dCNN}
\end{table}

\subsection{Two-dimensional CNN}\label{s.2dcnn}
We now apply a two-dimensional CNN to certain vector fields (visualised as digital images) determined by elliptic curves and Sato--Tate data.
To do so, we first apply the complex character twisting and image-generation pipeline from Section~\ref{s.twist} to the datasets $\mathcal{D}_1$ and $\mathcal{D}_2$. 
For $\mathcal{D}_1$, this pipeline yields a dataset $\widetilde{\mathcal{D}}_1$ of vector fields as in equation~\eqref{eq.vftci}.
For $\mathcal{D}_2$, it yields the analogous dataset $\widetilde{\mathcal{D}}_2$ in which $a_p\chi(p)$ is replaced by 
$\tilde x_{p} \chi(p)$.
The datasets $\widetilde{\mathcal{D}}_1$ and $\widetilde{\mathcal{D}}_2$ may be generated using code available at \cite{repo}. 
It follows from Section~\ref{s.1dcnn} that the datasets $\widetilde{\mathcal{D}}_1$ and $\widetilde{\mathcal{D}}_2$ are disjoint, and the combined dataset $\widetilde{\mathcal{D}}=\widetilde{\mathcal{D}}_1\sqcup\widetilde{\mathcal{D}}_2$ consists of 64,687 genuine elliptic curve twist arrays and 500,000 random matrix arrays.

By analogy with Section~\ref{s.1dcnn}, we set up the following  binary classification problem:
given $A\in\widetilde{\mathcal{D}}$, determine $i$ such that $A\in\widetilde{\mathcal{D}}_i$.
By construction, the images in $\widetilde{\mathcal{D}}_2$ match those in $\widetilde{\mathcal{D}}_1$ in pixel value distributions, image resolution, shape and normalisation.
They differ only in the arithmetic correlations between coefficients at different primes and different twists.
This setup ensures that any successful classification must rely on higher-order structure beyond Sato--Tate statistics, as captured by converse theorems.

We use a 2d CNN with architecture as described in Table~\ref{tab:CNN} and an 80:20 training ratio. 
We note that the architecture and hyperparameters are essentially the same as in Section~\ref{s.1dcnn}, with the obvious exception that the 1d convolutional filters have been changed to 2d. 
In words, after each convolutional layer, we apply batch normalisation, ReLU activation, and then max pooling to halve the spatial dimensions. After the final convolutional block, global average pooling reduces each feature map to a single scalar.
This is followed by a fully connected layer mapping
$512\mapsto 256$ features with batch normalisation, $\mathrm{ReLU}$ activation, and dropout (probability 0.5) for regularization.
A final linear layer produces a single logit for binary classification.
The network is trained using the Adam optimizer with a learning rate of $0.001$.
To address the class imbalance, we use the same conventions as in Section~\ref{s.1dcnn}.

\begin{table}[h]
\centering
\begin{tabular}{|c|c|c|}
\hline
Layer & Operation & Output shape \\
\hline
Input & -- & $(2, N, N)$ \\
Conv1 & $3\times3$ conv, 64 ch & \\
Pool1 & $2\times2$ max pool & $(64, [N/2], [N/2])$ \\
Conv2 & $3\times3$ conv, 128 ch & \\
Pool2 & $2\times2$ max pool & $(128, [N/4], [N/4])$ \\
Conv3 & $3\times3$ conv, 256 ch &  \\
Pool3 & $2\times2$ max pool & $(256, [N/8], [N/8])$ \\
Conv4 & $3\times3$ conv, 512 ch & \\
Pool4 & $2\times2$ max pool & $(512, [N/16], [N/16])$ \\
Conv5 & $3\times3$ conv, 512 ch & \\
Pool5 & $2\times2$ max pool & $(512, [N/32], [N/32])$ \\
GAP & Global average pooling & $(512)$ \\
FC1 & Fully connected & $(256)$ \\
FC2 & Fully connected & $(1)$ \\
\hline
\end{tabular}
\caption{CNN architecture used for classification problem described in Section~\ref{s.2dcnn}. Batch normalisation and $\mathrm{ReLU}$ activations follow each convolution.}
\label{tab:CNN}
\end{table}

\begin{figure}[htbp]
    \centering
    \includegraphics[width=0.6\textwidth]{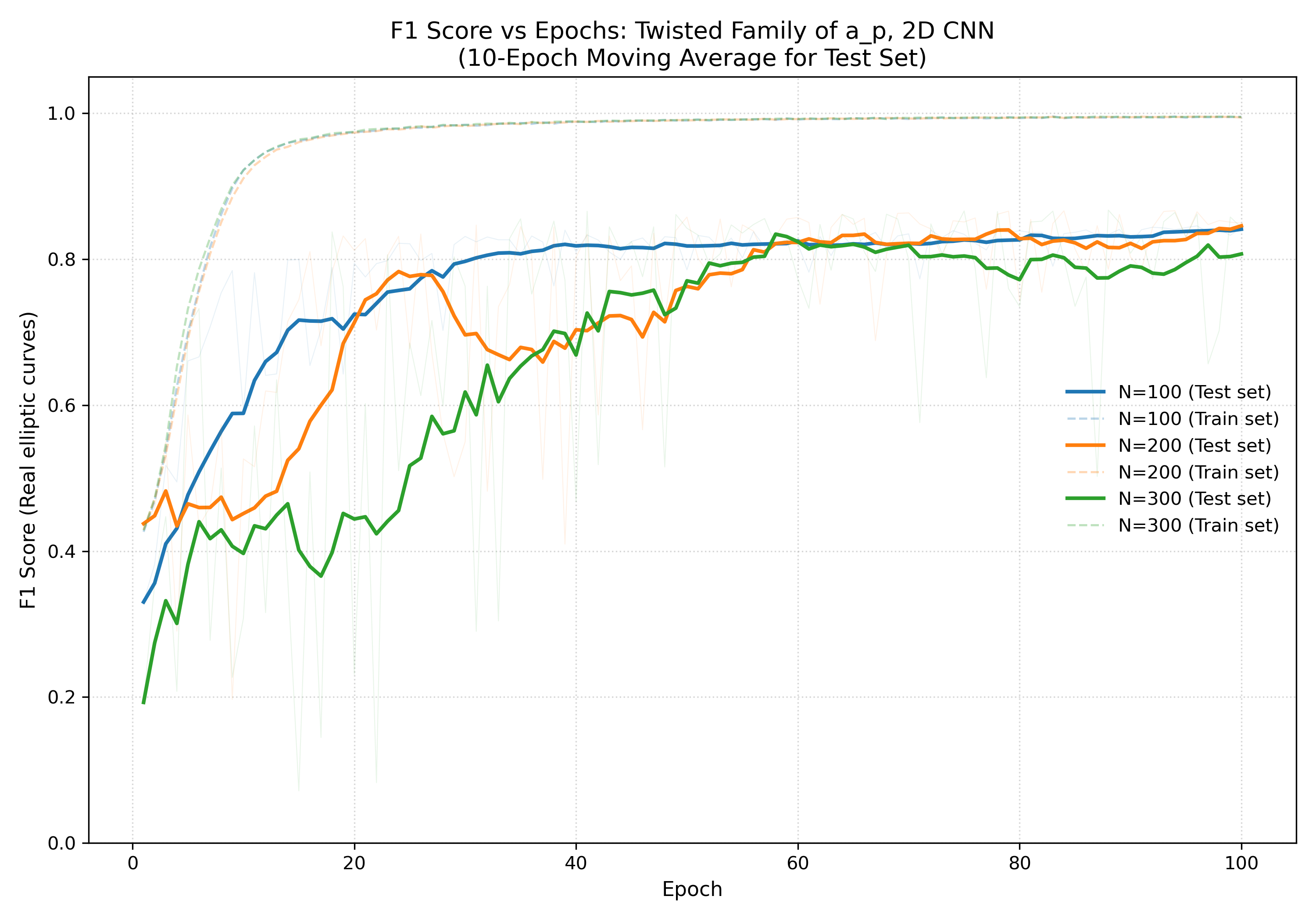}
    \caption{F1 scores on the training and test set for the classification problem in Section~\ref{s.2dcnn} using different numbers of Frobenius traces.}
     \label{fig:f1_100_200_300_comparison}
\end{figure}

For training, we use only the two non-constant channels (red and blue), and retain full numerical precision (rather than the visual encoding which involves rounding).
In Figure~\ref{fig:f1_1d_2d_comparison}, we see that the F1 score improves with the epoch, beginning at around $40\%$ and eventually reaching an F1 score exceeding $80\%$ after 100 epochs.
In Figure~\ref{fig:f1_100_200_300_comparison}, we see that the F1 score is similar for three different choices of $N\in\{100,200,300\}$.
It would be interesting to understand theoretically why the F1 score appears to stabilise around $80\%$. 
For example, it may be that the curves with larger conductor in our dataset simply require more twists to be verified as modular forms. 
This is a subtle question in the literature on converse theorems.
Indeed, for some conductors, it is known theoretically that no twists at all are needed \cite{CF,BCOZ}.
On the other hand, for arbitrary conductors, Weil's converse theorem requires an infinite family of twists by primitive Dirichlet characters \cite{weil}, whereas Razar's converse theorem requires only finitely many twists including imprimitive characters \cite{Razar}

\begin{figure}[H]
    \centering
        \includegraphics[width=0.3\textwidth]{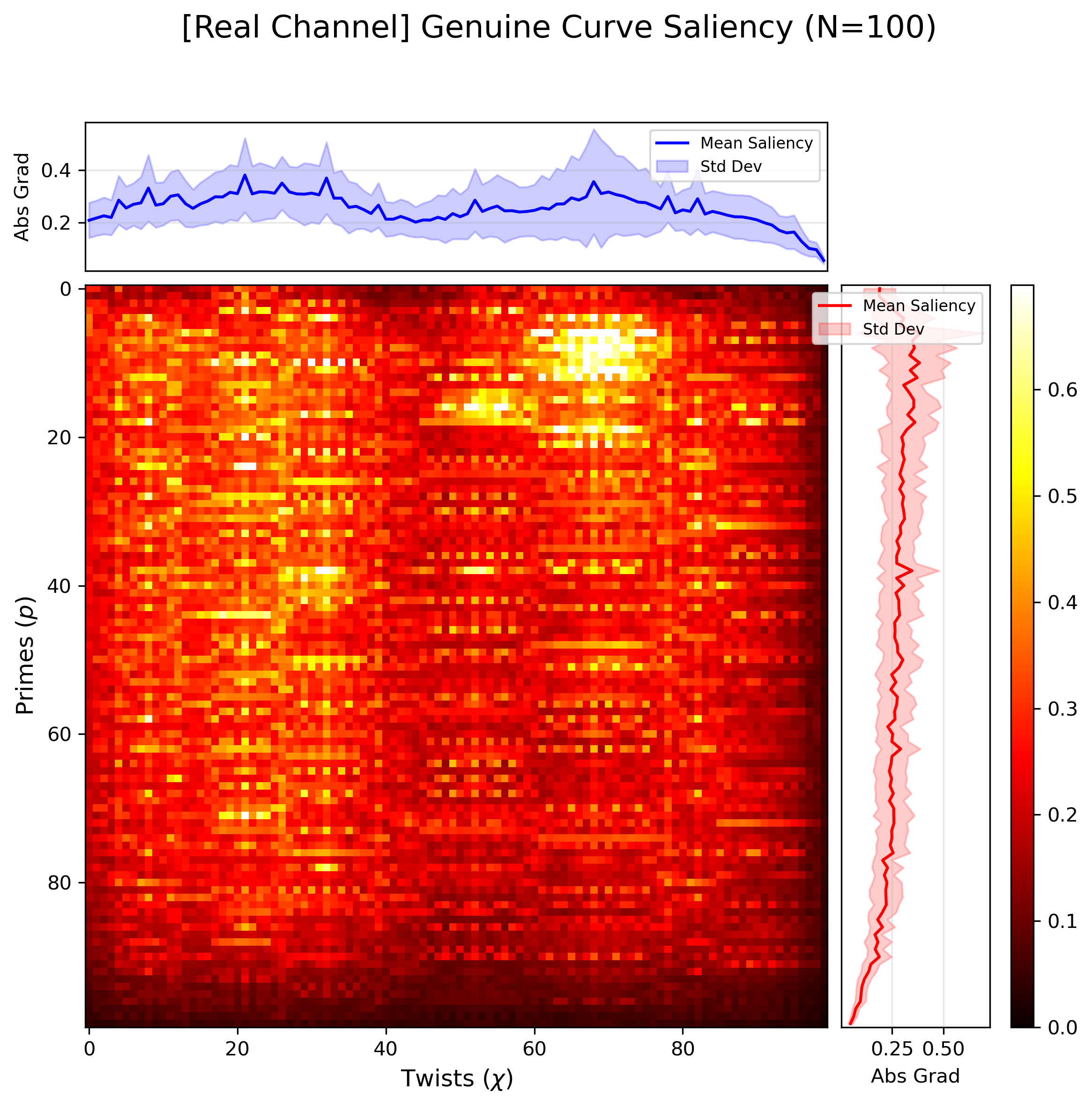}
        \includegraphics[width=0.3\textwidth]{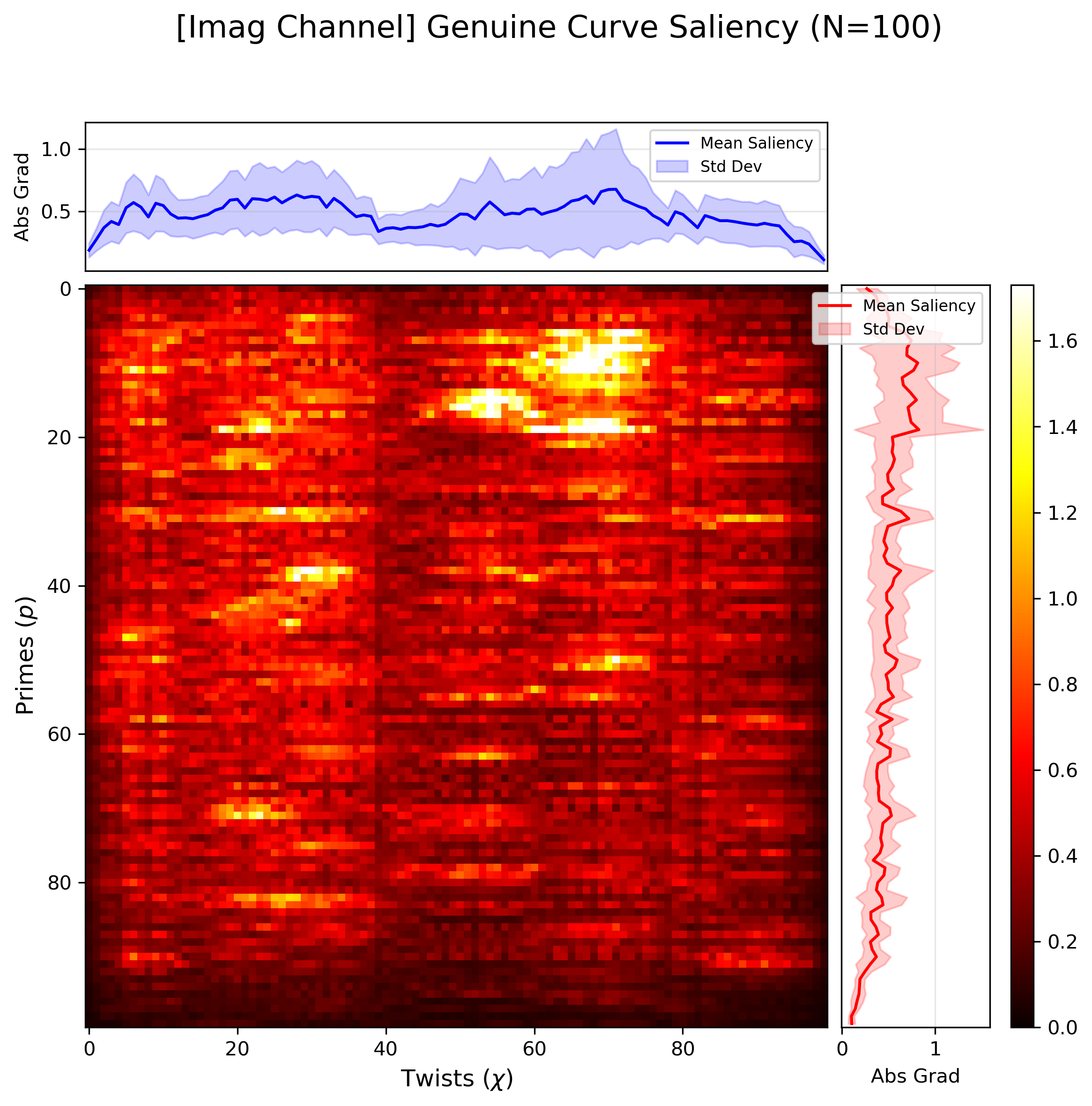}
        \includegraphics[width=0.3\textwidth]{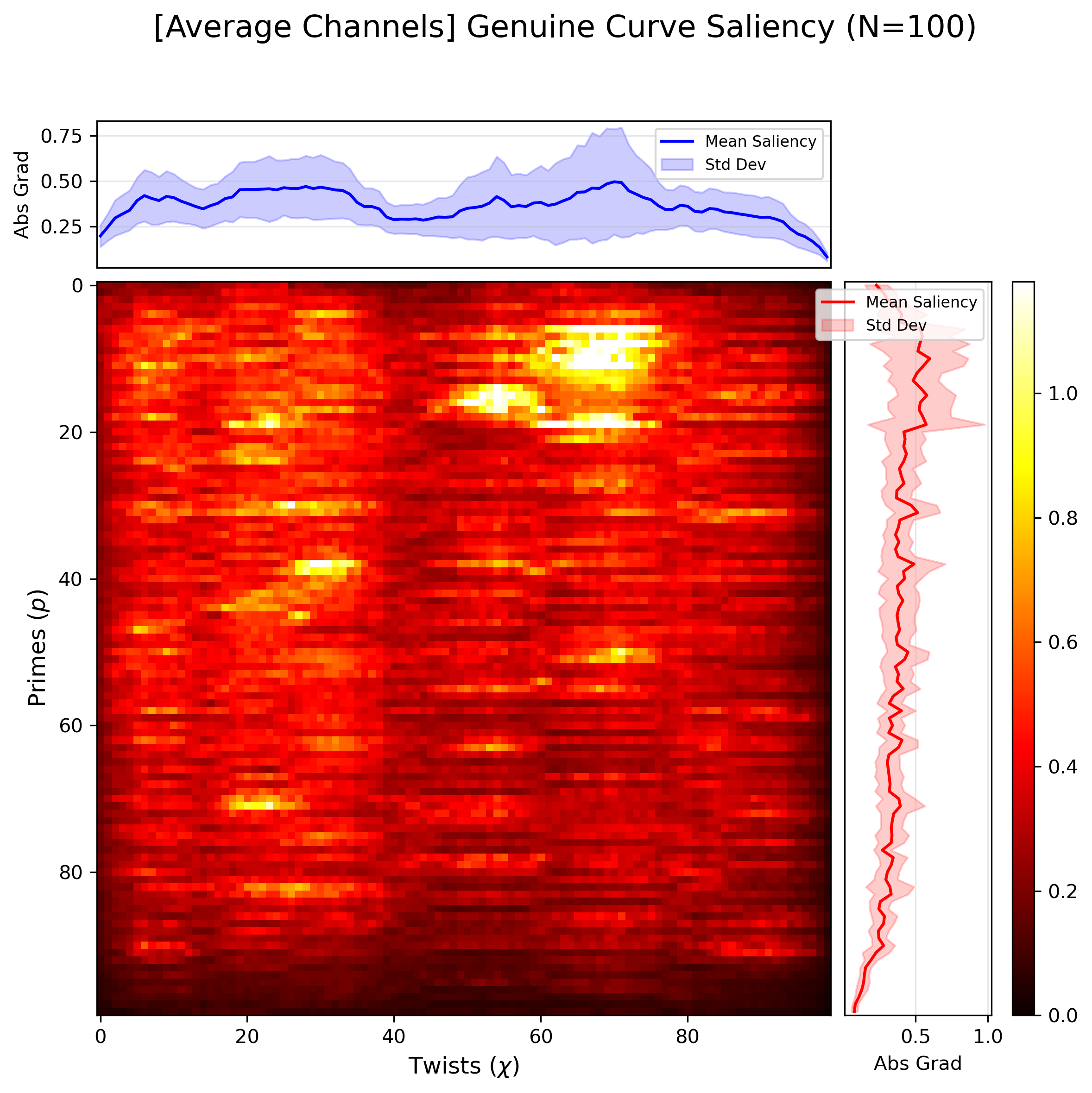}

        \includegraphics[width=0.3\textwidth]{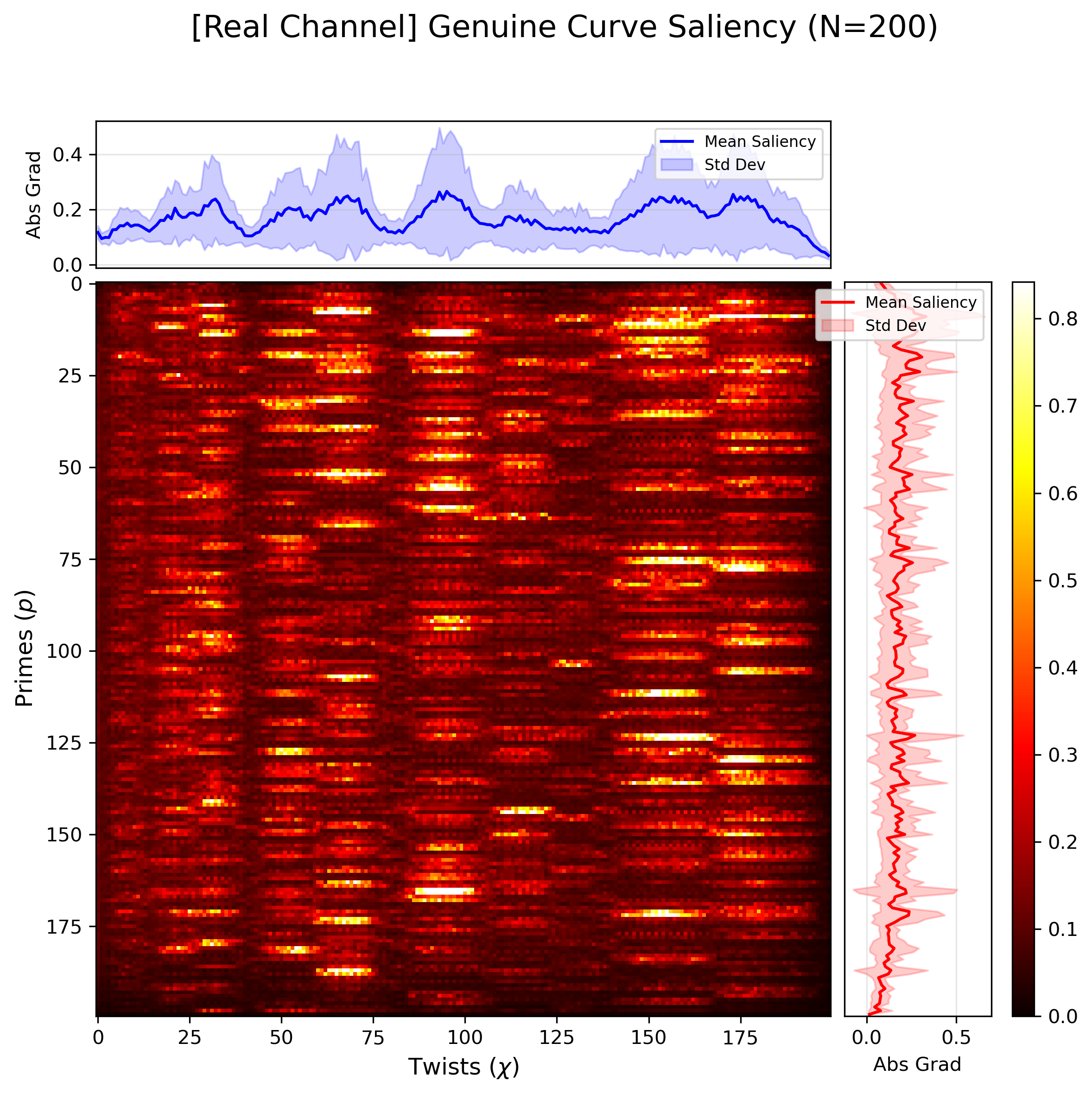}        
        \includegraphics[width=0.3\textwidth]{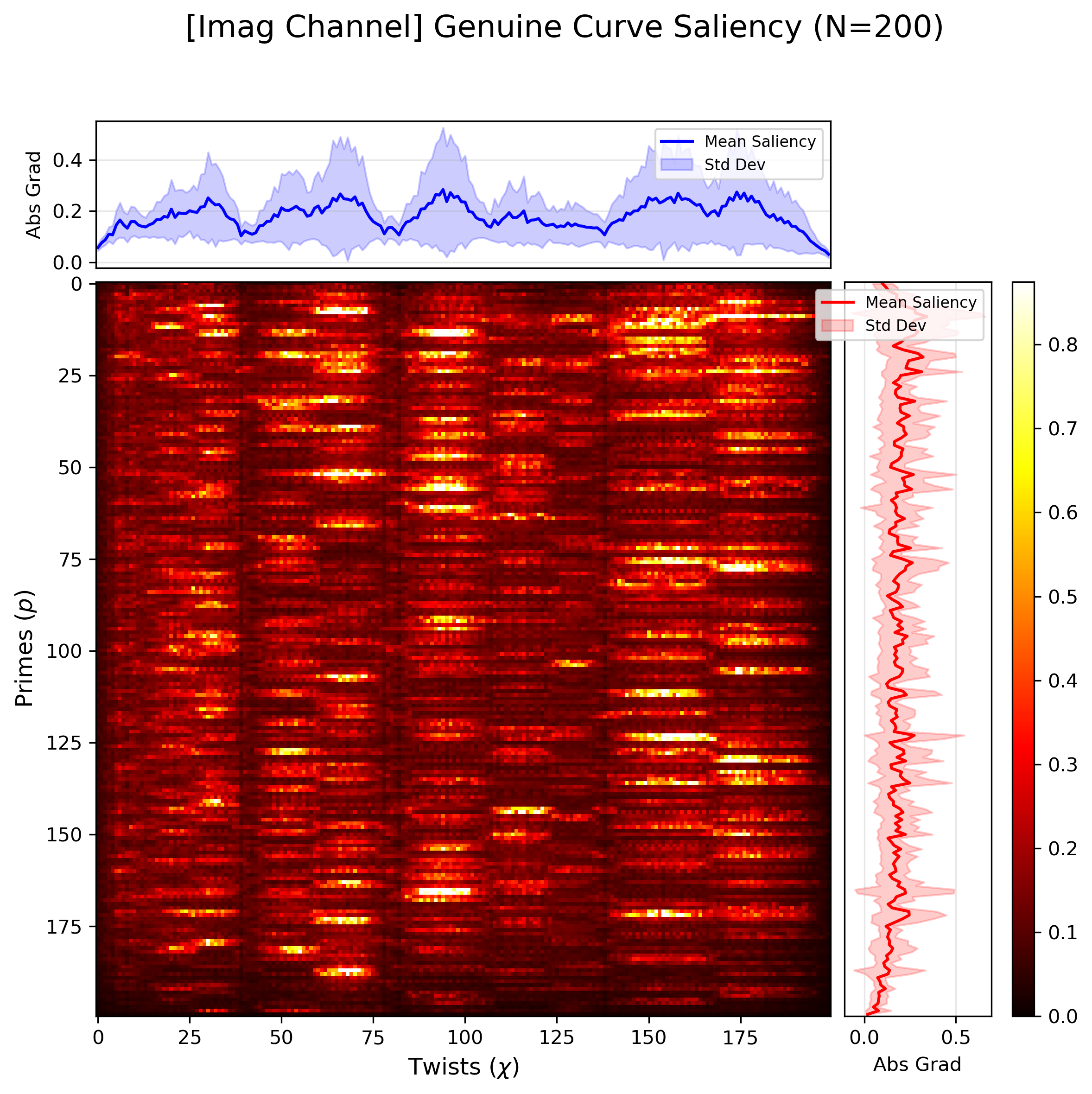}      
        \includegraphics[width=0.3\textwidth]{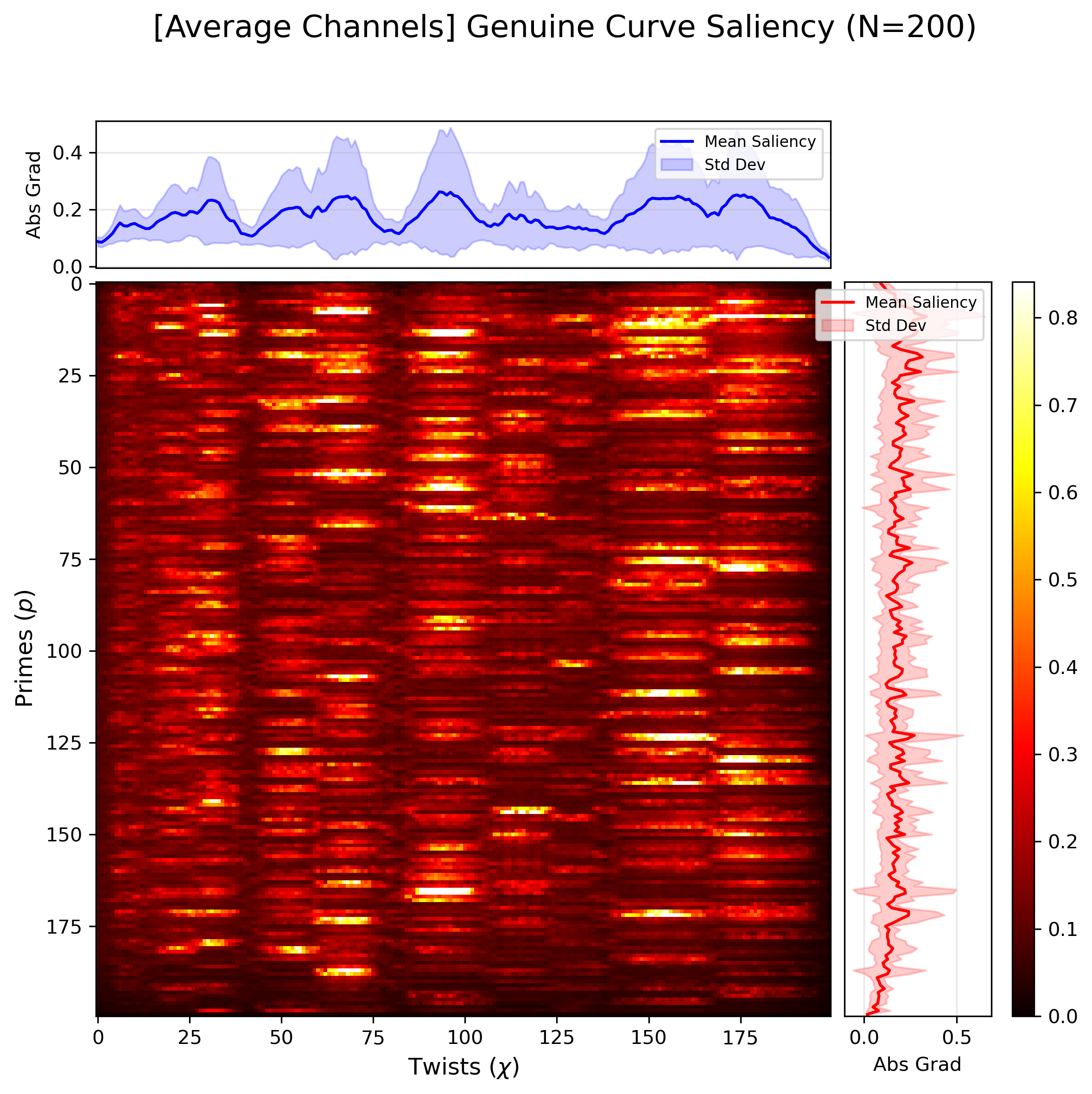}

        \includegraphics[width=0.3\textwidth]{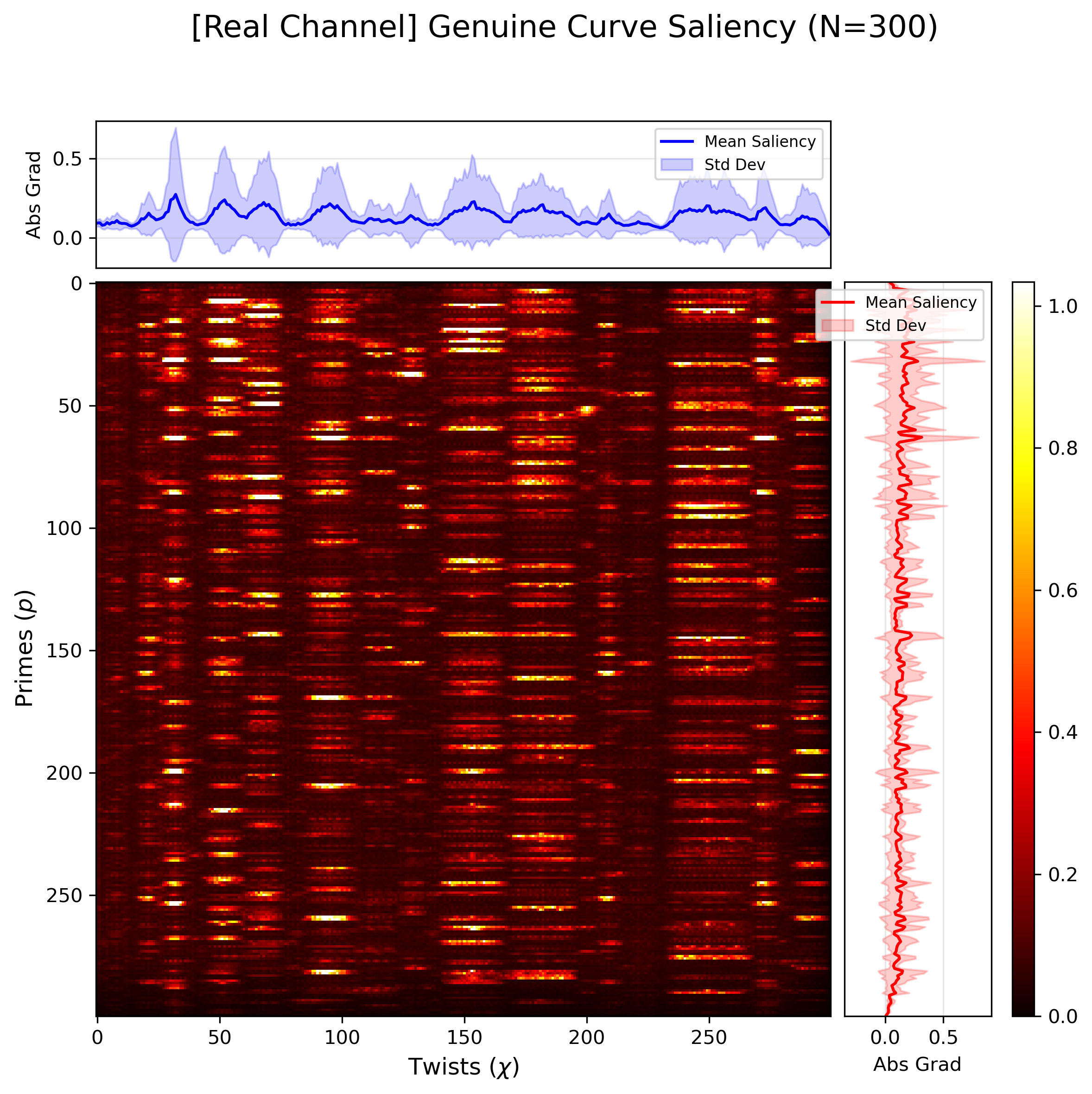}
        \includegraphics[width=0.3\textwidth]{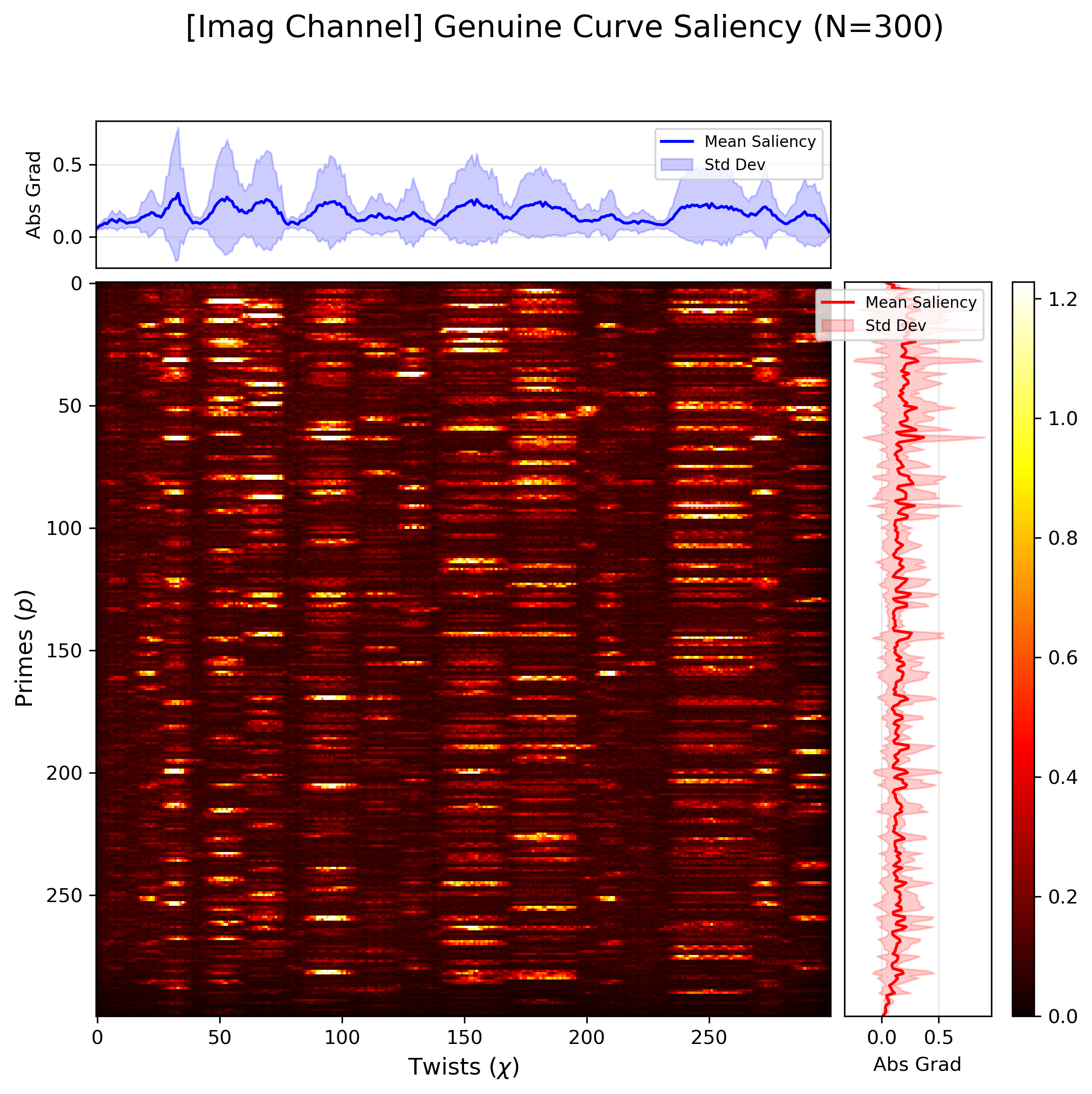}
        \includegraphics[width=0.3\textwidth]{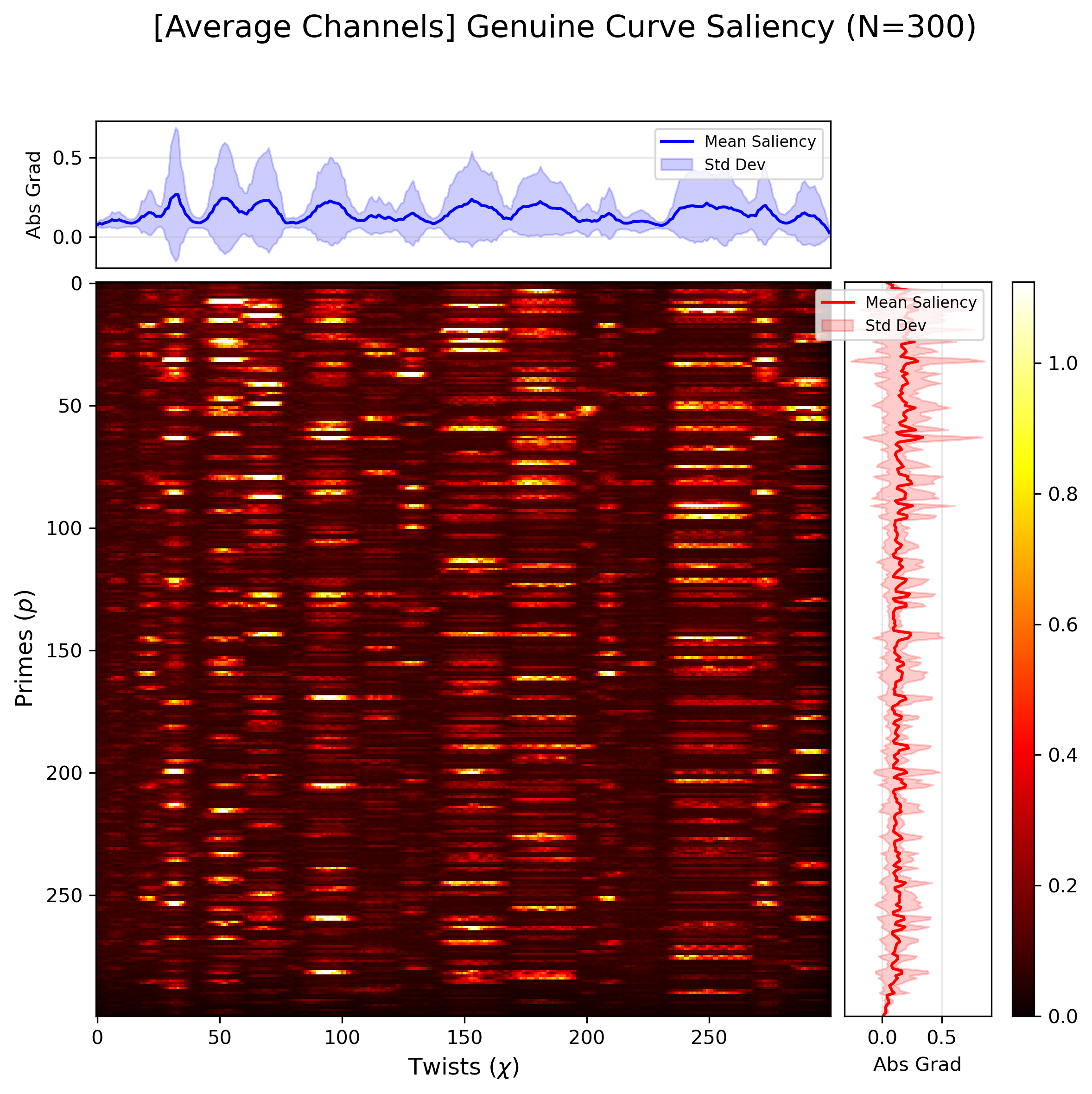}

    \caption{Two-dimensional saliency plots for the (left) real channel, (centre) imaginary channel, and (right) average of real and imaginary channels of the CNN described in Section~\ref{s.2dcnn}. In the top (resp. middle, bottom) row, the images have resolution $100\times100$ (resp. $200\times200$, $300\times300$).}
     \label{fig:saliency-real-imag-synthetic}
\end{figure}

In Figure~\ref{fig:saliency-real-imag-synthetic}, we exhibit (two-dimensional) saliency plots for our CNN, as we now explain.
The CNN determines a function
\[
c:S\rightarrow\{0,1\},
\] 
where $S$ is the set of $3$-dimensional vector fields on an $N\times N$ grid in $\mathbb{Z}^2$ as described in Section~\ref{s.twist}.
We may embed $S\hookrightarrow\mathbb{R}^{2N^2}$, in which the dimension is determined by the number ($N\times N$) of spatial parameters and the number ($2$) of intensity parameters.
We caution that $S$ has two (not three) intensity parameters because only two channels are non-constant by construction, and the intensity parameters take values in $[0,1]$ (not $\{0,\dots,255\}$) since we are working with the underlying vector field not the image encodings.
In fact, the CNN determines a sequence of such functions (one for each epoch), although we will suppress this from the notation.

The function $c$ factors through a function 
\[
f:S\rightarrow\mathbb{R},
\] 
which is composed with the sigmoid function and rounded to produce a value in $\{0,1\}$.
We define the saliency score at pixel $(p,\chi)$ in channel $R$ (resp. $B$) to be
\[S_R(p,\chi)=\frac{1}{|\widetilde{\mathcal{D}}_1|}\sum_{E\in\widetilde{\mathcal{D}}_1}\left|\frac{\partial f_E}{\partial R_{p,\chi}}\right|, \ \ \left(\text{resp. }S_B(p,\chi)=\frac{1}{|\widetilde{\mathcal{D}}_1|}\sum_{E\in\widetilde{\mathcal{D}}_1}\left|\frac{\partial f_E}{\partial B_{p,\chi}}\right|\right),\]
where $R_{p,\chi}$ (resp. $B_{p,\chi}$) are as in Section~\ref{s.twist}, $f_E$ denotes the evaluation of $f$, prior to the final sigmoid activation, at the image corresponding to $E$, and the partial derivatives are computed via automatic differentiation. 

In the heatmap component of Figure~\ref{fig:saliency-real-imag-synthetic}, we plot $S_R(p,\chi)$, $S_B(p,\chi)$ and $\frac12(S_R(p,\chi)+S_B(p,\chi))$, respectively.
The intensity of pixel in a saliency map indicates how much the corresponding pair $(\chi,p)$ affects the output of the CNN.
The blue curves indicate the average saliency across primes for each twist, and the red curve indicate the average saliency across twists for each prime.
The horizontal streaks in Figure~\ref{fig:saliency-real-imag-synthetic} indicate that the model focuses on some primes regardless of twist value, and the vertical clusters indicate that certain twists activate sensitivity across multiple primes.
The blue curves exhibit local minima near twist indices $39$ and $78$, which correspond to real Dirichlet characters, consistently across resolutions.
At such characters, $\chi(p)\in\{-1,0,1\}$, so the twisted data is merely a sign change of the untwisted Frobenius trace and the imaginary channel vanishes. We note that restricting to subsets of curves with prescribed rank  does not noticeably affect the saliency plots.
It would also be interesting to explore the consistency between the saliency vectors at different image resolutions.

\subsection{Higher conductor curves}\label{s.higherconductor}
\begin{figure}[H]
    \centering
    \includegraphics[width=0.8\linewidth]{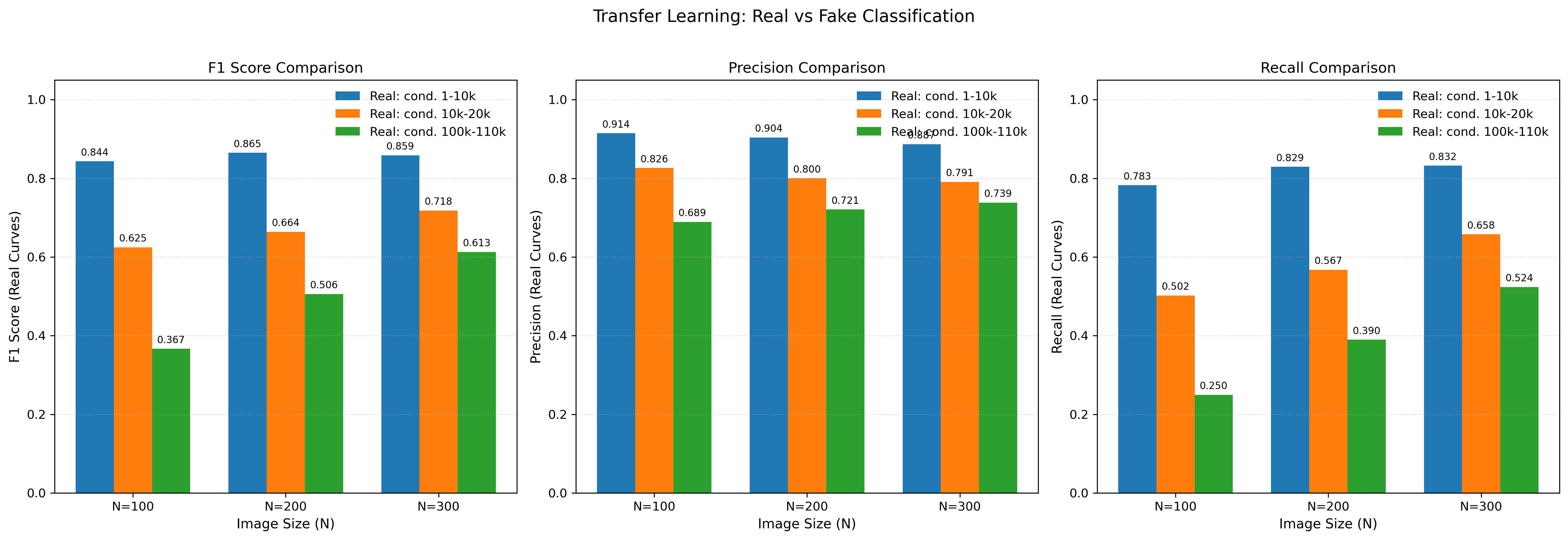}
    \caption{Evaluation metrics for transfer learning experiment described in Section~\ref{s.higherconductor}, namely (left) F1 score, (centre) precision, (right) recall.}
    \label{fig:transfer_f1}
\end{figure}

\begin{figure}[H]
    \centering
    \includegraphics[width=0.8\linewidth]{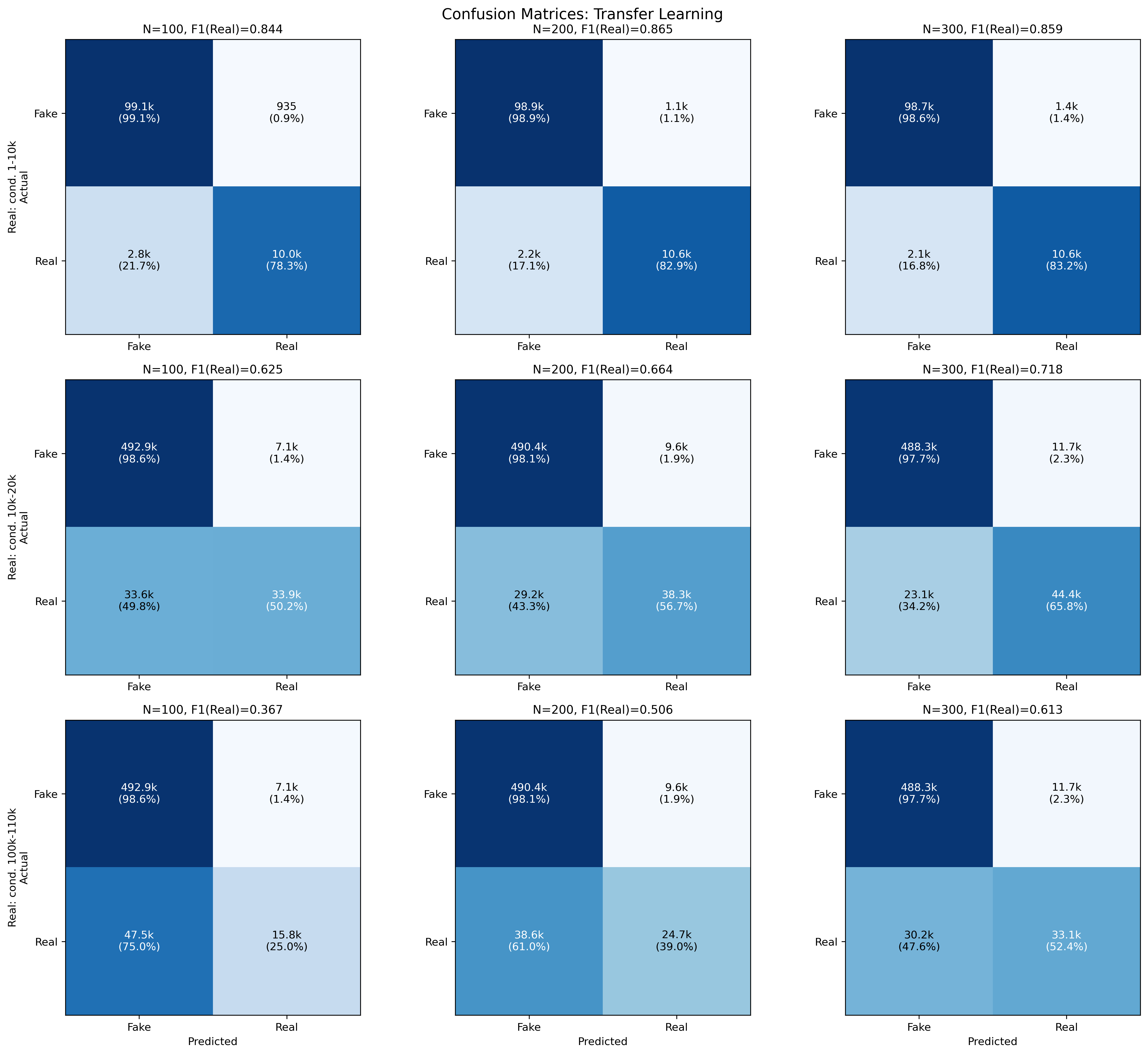}
    \caption{Confusion matrices for transfer learning experiment described in Section~\ref{s.higherconductor}}
    \label{fig:transfer_confusion_matrices}
\end{figure}

In Section~\ref{s.2dcnn}, we observed that the 2d CNN can distinguish curves with bounded conductor from random matrix data with the same distribution.
It is reasonable to ask if it is also capable of distinguishing unseen elliptic curves with larger conductor within a prescribed range.
To explore this question, we set up a transfer learning experiment in which we use the same architecture and training as in Section~\ref{s.2dcnn}, however, we evaluate the performance on entirely new test sets.
The test sets consist of new random matrix data generated as in Section~\ref{s.1dcnn} and genuine elliptic curves with conductor in the interval $[10000,20000]$ (resp. $[100000,110000]$).
In Figure~\ref{fig:transfer_f1}, we show the F1 score, precision, and recall for $N\in\{100,200,300\}$.
In Figure~\ref{fig:transfer_confusion_matrices}, we show the confusion matrices.

\subsection{Elliptic curve rank prediction}\label{s.rank}
\begin{figure}[H]
    \centering
    \includegraphics[width=0.47\textwidth]{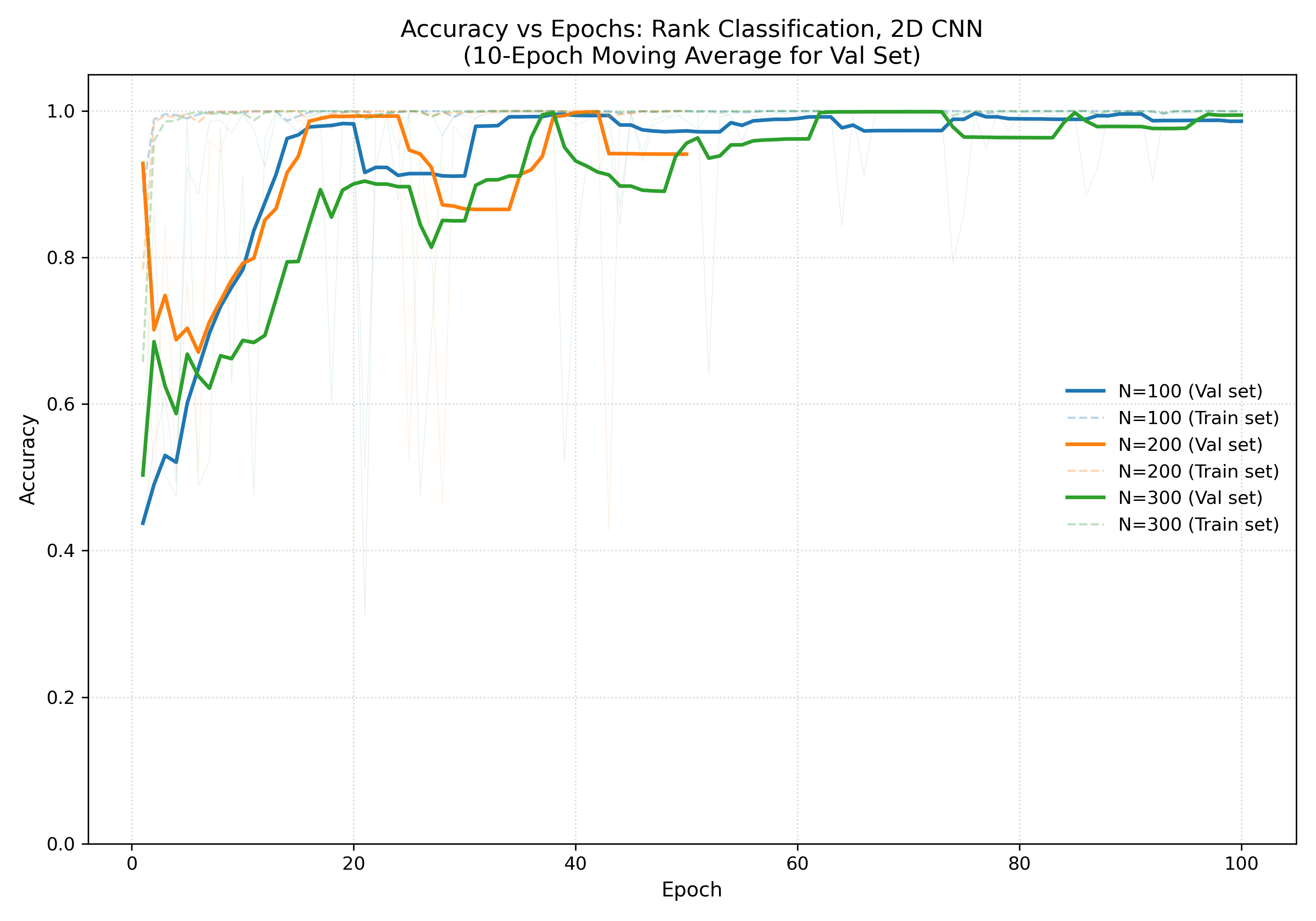}
    \caption{Evolution of accuracy for 2d CNN described in Section~\ref{s.rank} using different numbers of Frobenius traces.}
     \label{fig:accuracy-rank}
\end{figure}

\begin{figure}[H]
    \centering
    \includegraphics[width=0.47\textwidth]{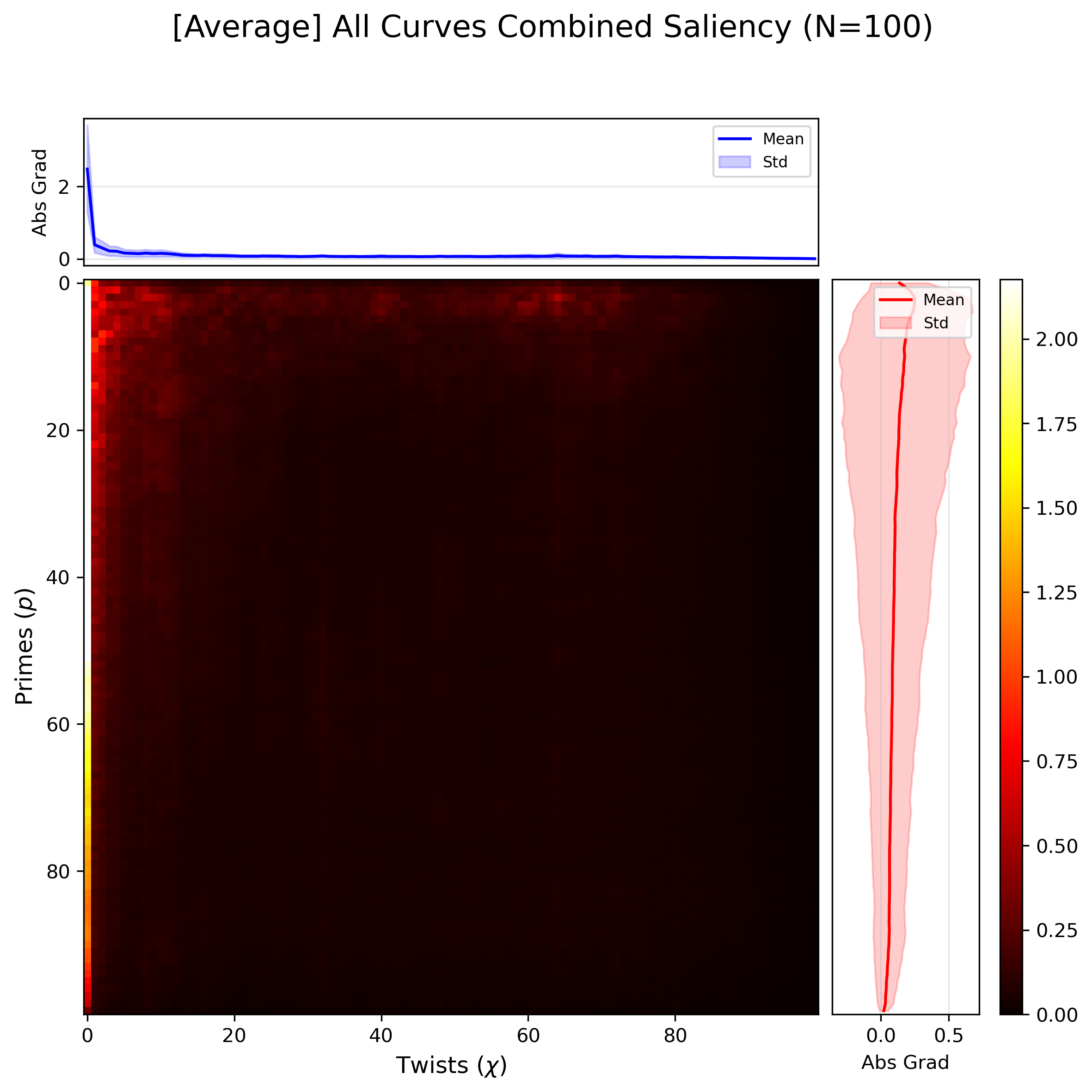}
    \caption{Two-dimensional saliency map for 2d CNN described in Section~\ref{s.rank}. }
     \label{fig:saliency-rank}
\end{figure}

A recent trend in the literature is to machine learn the rank of elliptic curves from Frobenius traces.
Whilst preliminary investigations used classical machine learning \cite{Alessandretti:2019jbs, He:2020tkg}, CNNs were applied to this task in \cite{KV22, Poz, big-group}.
In this section, we seek to machine learn a classification
\[R:\widetilde{\mathcal{D}}_1\rightarrow\mathbb{Z}_{\geq0},\]
where $\widetilde{\mathcal{D}_1}$ is as in Section~\ref{s.2dcnn}. 
The classification $R$ is intended to approximate the map $r:\mathcal{E}\rightarrow\mathbb{Z}_{\geq0}$, in which each elliptic curve $E\in\mathcal{E}$ is assigned to its rank $r(E)$.
This dataset $\widetilde{\mathcal{D}_1}$ has never been used before in the literature, and the main difference between this classification task and those in previous studies is the nature of the input data.
We use an almost identical 2d CNN architecture as in Section~\ref{s.2dcnn} with the only difference being that the final layer now produces outputs in $\{0,1,2\}$ according to the predicted rank.
In Figure~\ref{fig:accuracy-rank}, we show the evolution of the accuracy of the model on the validation set for $N\in\{100,200,300\}$.
We observe that the accuracy approaches 100\%.
In Figure~\ref{fig:saliency-rank}, we plot the 2d saliency maps associated to this classification. 
We observe a bright vertical line along the first twist, indicating that the CNN understands the redundancy of the information and learns to use only one column to predict the rank.

\bibliography{references}
\bibliographystyle{unsrt}

\end{document}